\documentclass[12pt]{article}
\usepackage[utf8]{inputenc}
\usepackage{amsthm}
\usepackage{amsmath}
\usepackage{amssymb}

\def\A{\mathcal A}
\def\B{\mathbb B}
\def\C{\mathbb C}
\def\H{\mathbb H}
\def\har{\mathcal H}
\def\M{\mathcal M}
\def\N{\mathcal N}
\def\R{\mathbb R}
\def\im{\,\!\mbox{Im}\,}
\def\re{\mbox{Re}\,}
\def\sc{\mbox{Sc}\,}
\def\Vec{\mbox{Vec}\,}
\let\cg\overline 
\def\Dbar{\cg D}
\newcommand{\scpr}[2]{\langle#1,\,#2\rangle}
\def\iff{\Longleftrightarrow}
\def\ker{\mbox{\,\!ker}\,}

\def\DS{\displaystyle}
\def\str{\raisebox{1.5ex}{\rule{0pt}{1ex}}}

\def\qed{\hspace{1em}\raise3.5pt\hbox{\framebox[1.3ex]{\rule{0pt}{.7ex}}}\vspace{1.5ex}}

\numberwithin{equation}{section}
\theoremstyle{plain}
\newtheorem{theo}{Theorem}[section]
\newtheorem {lemm}[theo]{Lemma}
\newtheorem {prop}[theo]{Proposition}
\newtheorem {coro}[theo]{Corollary}
\def\defi{\refstepcounter{theo}\par\noindent
 \textbf{Definition \arabic{section}.\arabic{theo}. }}
\def\proof{\par\vspace{-1ex}\noindent\textit{Proof. }}
\def\note{\refstepcounter{theo}\par\noindent
 \textbf{Note \arabic{section}.\arabic{theo}. }}

\begin{document}

\begin{center}{\Large Contragenic Functions of Three Variables }
  \\[2ex] Cynthia Álvarez-Peña$\hspace{-1pt}{}^1$  \\ 
Depto. de Matem\'aticas, Cinvestav-I.P.N.\\ 
  Quer\'etaro, M\'exico \\[1ex]
   R. Michael Porter\footnote{Partially supported by
CONACyT grant 166183}\\ 
Depto. de Matem\'aticas, Cinvestav-I.P.N.\\ 
  Quer\'etaro, M\'exico \\[1ex]
\end{center}

\noindent\textbf{Abstract:} It is shown that harmonic functions from a
simply connected domain in $\R^3$ to $\R^3$ cannot always be
expressed as a sum of a monogenic (hyperholomorphic) function and an
antimonogenic function, in contrast to the situation for complex
numbers or quaternions.  Harmonic functions orthogonal in
$L_2$ to all such sums are termed ``contragenic'' and their
properties are studied. A ``Bergman kernel'' and is derived, whose
corresponding operator vanishes precisely on the contragenic
functions.  A graded orthonormal basis for the contragenic function in
the ball $\B^3$ is given.



\setcounter{section}{-1}
\section{Introduction}

The following fact in elementary complex analysis is well known.

\textbf{Theorem A.} \textit{ Every harmonic function $u\colon
  \{|z|<1\} \to\C$ of a complex variable is expressible as the sum of
  a holomorphic function and an antiholomorphic function.}

This principle has many uses.  In particular, when a holomorphic
solution to a problem is sought and a first attempt is made as a
harmonic function, one may ``throw away'' the antiholomorphic part to
obtain a holomorphic approximation.  A classical example of this
principle in conformal mapping theory is found in the method of
Fornberg \cite{Fo}, in which a guess as to the boundary values of the
mapping is expressed on the boundary of the unit disk as a Fourier
series, whose coefficients give the sum of a power series in $z$
(positive powers of $e^{i\theta}$) and in $\cg z$ (negative powers) in
the interior of the disk. The sum of these two is a harmonic
function, whose antiholomorphic part is discarded in the algorithm.
 
Theorem A holds in many generalizations of the field $\C$ of complex
numbers, for example monogenic (hyperholomorphic) functions on
quaternions \cite{Sud} or on Clifford algebras \cite{BDS,GS}.  It
also holds for monogenic functions from $\R^3$ to $\H$
\cite{Ca,GHS}.

In this paper we show (Theorem \ref{theo:contragenic}) that the
natural generalization of Theorem A does not hold for monogenic
functions from $\R^3$ to $\R^3$.  Therefore, the class of harmonic
functions which have no such decomposition is new. We will use the
term ``contragenic'' for harmonic functions which are orthogonal to
the monogenic and the antimonogenic functions in the sense of
$L_2$. We believe they are of interest because of the relevance of
monogenic functions in $\R^3$ to physical systems \cite{SSS}.

Thus we initiate here the study of contragenic functions. The precise
definitions for monogenic functions and related notions in $\R^3$ are
specified in section \ref{sec:monofunc}, based on work of
\cite{Bock2,Ca,MG,Mor}. In section \ref{sec:homomono} we summarize the
necessary facts about the standard basis for homogeneous monogenic
polynomials in the unit ball of $\R^3$ and related spaces. We also
calculate an orthonormal basis for the vector parts of the monogenic
functions. In section \ref{sec:contragenic} we prove the existence of
contragenic functions and prove some of their basic properties. In
particular we derive a Bergman kernel which annihilates precisely the
contragenics.  Finally, in section \ref{sec:bases} we give an explicit
construction for a graded basis for the space of contragenic
functions.

\section{Monogenic functions in $\R^3$\label{sec:monofunc}}

\subsection{Notation}
We will use fairly standard notation for the skew field of quaternions
$\H=\R^4=\{x=x_0e_0+ x_1e_1+x_2e_2+x_3e_3,\ x_j\in\R\}$. Here $e_0=1$
is the unit and the multiplication is determined by $e_{j_1}^2=-1$ and
$e_{j_1}e_{j_2}=\pm e_{j_3}$ where $\{j_1,j_2,j_3\}$ is the set of
indices $\{1,2,3\}$ and the $+$ sign is taken precisely when the
cyclic order matches that of 1,2,3. See texts such as \cite{GHS,K,Sud}
for further general information.  The $\R$-subspace of quaternions
$x\in\H$ such that $e_3=0$ will be denoted $\R^3\oplus\{0\}$ or more
simply just $\R^3$ when there is no danger of confusion.

Let us write $\partial_j=\partial/\partial x_j$ for $j=0,1,2,3$.
There is a great deal of literature (see \cite{Sud} and references in
\cite{BDS}) on the Cauchy-Riemann type differential operators
operators
\begin{eqnarray*} \vec\partial_3 &=& \sum_{j=1}^3 \partial_j e_j,  \\
  D_\H     &=&  \partial_0 - \vec\partial_3, \\
  \Dbar_\H &=&  \partial_0 + \vec\partial_3.
\end{eqnarray*}
which act both from the left and from the right on differentiable
functions $f=f_0 +f_1e_1+f_2e_2+f_3e_3$ defined in open subsets of
$\H$.  The operator $\Dbar_\H$ is a generalization of the operator
$\partial/\partial \cg z$ on which complex analysis is based:
functions for which $\Dbar_\H f=0$ (resp.\ $f\Dbar_\H=0$) are variously
called left (resp.\ right) Fueter-regular, monogenic, hyperholomorphic, among others;
occasionally the roles of $D_\H$ and $\Dbar_\H$ are interchanged in
the terminology.

In recent years some work has been done \cite{Bock2,GHS} on the
analogous functions from $\R^3$ to $\R^4$ and from $\R^4$ to $\R^3$,
with a view to expressing and studying operators relevant to physics.
Relatively little has been done for functions from $\R^3$ to $\R^3$;
in this regard we mention \cite{GLS,MG}.  This setting is particularly
interesting because even though such functions are not conserved under
multiplication by elements of $\R^3\oplus\{0\}$, i.e., the algebraic
structure of an algebra or ring is lost, in many ways the monogenic
functions behave more like standard holomorphic functions in $\C$
(cf.\ Proposition \ref{prop:leftright}). To be precise, for $x=x_0 +
x_1e_1+x_2e_2\in\R^3$, let us write $\sc x=x_0$, $\Vec x =\vec x =
x_1e_1+x_2e_2$, $\cg x= x-\vec x$. Let $\Omega\subset\R^3$ be an open
set, and $\har_\R(\Omega)$ the space of real-valued harmonic functions
defined in $\Omega$. We consider the (Moisil-Theodorescu type) operators $D$
and $\Dbar$, defined by
\begin{equation}
 \vec\partial =  \partial_1 e_1+ \partial_2 e_2, \quad
  D = \partial_0 - \vec\partial, \quad
  \Dbar = \partial_0 + \vec\partial ,
\end{equation}
and define the set of \textit{(left-)monogenic} (or
 \textit{hyperholomorphic}) functions
\begin{equation}
 \M(\Omega)= \{f= f_0 +f_1e_1+f_2e_2 \in C^1(\Omega,\R^3)\colon\ \Dbar f=0\}.
\end{equation}
The fact that $\M(\Omega)\subseteq\har(\Omega)$, where
$ \har(\Omega)=\har_\R(\Omega)\times \har_\R(\Omega)\times \har_\R(\Omega)$ is the set of $\R^3$-valued harmonic functions in $\Omega$, follows immediately from the
factorization $\Delta = D\Dbar=\Dbar D$ of the Laplacian on $\R^3$. 
 
\subsection{Antimonogenics and ambigenics\label{subsec:antiambi}}
 We say that $f$ is \textit{(left) antimonogenic} when $Df=0$. 
  A \textit{monogenic constant}  is a function which is simultaneously
monogenic and antimonogenic:  $Df=\cg Df=0$, or equivalently,
$\partial_0f=\vec\partial f=0$.

It is unavoidable that $Df$ and $\cg Df$ need not take their values in
$\R^3$ even when $f$ does. However, due to the fact that
\begin{eqnarray*}
 -e_3 \cg D f e_3  &=&
 -e_3(\partial_0+\vec\partial) e_3(-e_3)(f_0+\vec f)e_3  \\
 &=& (\partial_0-\vec\partial)(f_0-\vec f)\\
 &=& D \cg f
\end{eqnarray*}
we have
\[    \cg D f =0 \iff D \cg f=0 \iff f \cg D=0,
\]
and consequently,
\begin{prop} \label{prop:leftright}
A function is left monogenic if and only if it is right monogenic. The
set of conjugates of monogenic functions
\begin{equation}  \cg\M(\Omega)= \{ \cg f\colon\ f\in\M(\Omega) \}   
\end{equation} 
coincides with the set of antimonogenic functions in $\Omega$.
\end{prop}
 
This is of course quite different from the situation for monogenic
functions in $\H$, where left- and right-monogenicity are
different. The proposition allows us to write
$\M(\Omega)\cap\cg\M(\Omega)$ for the set of monogenic constants.  If
$f\in\M(\Omega)\cap\cg\M(\Omega)$, then $\partial_0f=\vec\partial
f=0$. Thus a monogenic constant $f$ does not depend on $x_0$ and can
be expressed as
\begin{equation} \label{eq:constmono}
 f= c_0+f_1e_1+f_2e_2
\end{equation}
where $c_0\in\R$ is constant and the quantity $f_1-if_2$ is an
ordinary holomorphic function of the complex variable $x_1+ix_2$. 
There are natural projections of $\M(\Omega)$ onto the subspaces
\begin{eqnarray*}
 \sc\M(\Omega) &=& \{\sc f\colon\ f\in\M(\Omega)\}
 \ \subseteq\ \har_\R(\Omega),\\
 \Vec\M(\Omega) &=& \{\Vec f\colon\ f\in\M(\Omega)\}
 \ \subseteq\ \har_{\{0\}\oplus\R^2}(\Omega),
\end{eqnarray*}
and by Proposition \ref{prop:leftright} we see that
$\sc\M(\Omega)=\sc\cg\M(\Omega)$, $\Vec\M(\Omega)=\Vec\cg\M(\Omega)$.
An element of $\M(\Omega)+\cg \M(\Omega)$ will be called an
\textit{ambigenic} function; its decomposition as a sum of a
monogenic and an antimonogenic function is unique up to the addition
of a monogenic constant.

Consider $L_2(\Omega,\R^3)$ with the real-valued inner product
$\scpr{f}{g}=\sc\int_{\Omega}\cg fg\,dV
=\int_{\Omega}(f_0g_0+f_1g_1+f_2g_2)\,dV$. Write
\[ \M_2(\Omega)= \M(\Omega)\cap L_2(\Omega).
\]
It is somewhat inconvenient that $\cg\M_2(\Omega)$ is not orthogonal
to $\M_2(\Omega)$.  However, we have automatically that
$\sc\M_2(\Omega)\perp \Vec\M_2(\Omega)$ since any scalar function
multiplied by $e_0$ is by definition orthogonal to any combination of
$e_1$ and $e_2$. This fact gives us an orthogonal direct sum
decomposition of the space of square-integrable ambigenic functions,
\begin{equation} \label{eq:dirsum}
 \M_2(\Omega)+\cg \M_2(\Omega) = \sc\M_2(\Omega) \oplus \Vec\M_2(\Omega).
\end{equation}

We will always assume that $\Omega\subseteq\R^3$ is connected.
\begin{lemm} \cite{Mor} \label{lemm:f0} 
Suppose $\Omega$ is simply connected. 
(a) Let
$f_0\in \har_\R(\Omega)$. Then there exists $f\in\M(\Omega)$ such that
$\sc f=f_0$. This $f$ is unique up to an additive monogenic constant.

(b) Let $f_1,f_2\in \har_\R(\Omega)$. A necessary and sufficient
 condition for there to exist $f=f_0+\vec f\in\M(\Omega)$ such that
 $\sc f=f_0$, $\vec f=f_1e_1+f_2e_2$, is that $\partial_2
 f_1=\partial_1 f_2$.  When it exists, this $f$ is unique up to an
 additive scalar constant.
\end{lemm}

When the operator $D$ is applied exclusively to scalar-valued harmonic functions,
i.e.
\[    D \colon \har_\R(\Omega) \to  \M,
\]
we see from (\ref{eq:dirsum}) and Lemma \ref{lemm:f0} that this
operator splits naturally to give exact sequences
\begin{equation} \label{eq:par0}
    0\to  \R \to \har_\R(\Omega) \stackrel{\partial_0}{\to}  \sc\M\to 0, \\
\end{equation}
\begin{equation}\label{eq:parvec}
  0\to \ker\vec\partial\to \har_\R(\Omega) \stackrel{ \vec\partial}{\to}  \Vec\M \to0. 
\end{equation}
Here $\ker\vec\partial$ is two-dimensional over $\R$, consisting only of
polynomials $h(x)=c_0+c_1x_0$.

\note One could equally well embed $\R^3$ in $\H$ differently, for
example by considering $ \{0\} \oplus \R^3 \subseteq\H$, thus writing
\[
\begin{array}{ccccrcr}
 \tilde x &\!=\!& \tilde x_1e_1+\tilde x_2e_2+\tilde x_3e_3 
 & \!=\!&  -e_3x  &\!=\!&   x_2e_1 - x_1e_2+x_0e_3   \\
 \tilde f &\!=\!& \tilde f_1e_1+\tilde f_2e_2+\tilde f_3e_3 
 &\!=\!& \ \; e_3f  &\!=\!& -f_2e_1 + f_1e_2 -f_0e_3  
\end{array}
\]
and using the operators 
\[ \widetilde D = 
   e_1 \tilde\partial_1-e_2 \tilde\partial_2-e_3 \tilde\partial_3 ,
 \quad  \cg{\widetilde D} = 
   e_1 \tilde\partial_1+e_2 \tilde\partial_2+e_3 \tilde\partial_3.
\]
The relationship $ \cg{\widetilde D}\tilde f |_{\tilde x}  =  Df|_{e_3\tilde x} $
implies
 \begin{eqnarray*}
    D f =0  &\iff&  \cg{\widetilde D} \tilde f =0 \\
  &\iff&    \mbox{grad}\,\tilde f =0 ,\ \mbox{curl}\,\tilde f =0 
\end{eqnarray*}
and consequently this alternative embedding, which is used for
example in \cite{GS,GLS}, is equivalent to the form we are using in
this article. Our spaces $\M(\Omega),\cg\M(\Omega)$ correspond to the
spaces of left- and right-monogenic functions $\tilde f$ in the sense
of $\cg{\widetilde D}$. In this context the equations defining
monogenicity are also known as a Riesz system \cite{De,Mor,MG,Rie}.

\section{Homogeneous monogenics\label{sec:homomono}}
  
We will mostly work in the ball $\B^3=\{|x|<1\}\subseteq\R^3$. The
real vector space $\har_\R(\B^3)\cap L_2(\B^3)$ is stratified into the
subsets $\har^{(n)}_\R(\B^3)$ comprised of real harmonic
functions which are homogeneous of successive degrees
$n=0,1,\dots$\ It is well known that the elements of $\har^{(n)}_\R(\B^3)$ 
are polynomials of degree $n$ and that these real linear
subspaces are orthogonal with respect to the inner product
$\langle\cdot,\cdot\rangle$.  Further, every square integrable
harmonic function on $\B^3$ has a unique expression as a series formed of
elements of these sets.

\subsection{Basis for $\M_2(\B^3)$}
Many schemes have been devised to construct bases for spaces of
homogeneous monogenic functions of given degree, from the classical
construction of Fueter to diverse applications of symmetric sums of
products; see \cite{Ca,GHS,Ma}. For $\R^3$ one way is to proceed as
follows. A well known orthonormal basis of $\har^{(n)}_\R(\B^3)$,
$n\ge0$, is the system of $2n+1$ solid spherical harmonics
\begin{equation}
   \widehat U_0^n, \widehat U_1^n,\dots, \widehat U_n^n,\widehat V_1^n,\dots,\widehat V_n^n,
\end{equation}
where $\widehat{U}^{n}_m= r^n U^{n}_m$, $\widehat{V}^{n}_m= r^n
V^{n}_m$, are defined in terms of spherical coordinates $x_0=r \cos
\theta$, $x_1=r \sin \theta \cos\varphi,$ $x_2=r\sin\theta
\sin\varphi$ via the relations
\begin{eqnarray} \label{eq:defUV}
 U_0^n(\theta,\varphi) &=& P_n(\cos\theta), \nonumber\\
 U_m^n(\theta,\varphi) &=& P_n^m(\cos\theta)\cos(m\varphi), \nonumber\\
 V_m^n(\theta,\varphi) &=& P_n^m(\cos\theta)\sin(m\varphi), \qquad m=1,2,\cdots,n.
\end{eqnarray}
Here $P_n$ is the Legendre polynomial of degree $n$ and $P_n^m$ is
the associated Legendre function is given by
\[   P_n^m(t)=(1-t^2)^{m/2} \frac{d^m}{dt^m}P_n(t).
\] 
We will need the following standard identities \cite{Olv}:
\begin{eqnarray}
  (1-t^2)(P^m_{n+1})'(t)  &=& (n+m+1) P^m_n(t) - (n+1)t P^m_{n+1}(t) ,  \label{eq:L1} \\
  (1 - t^2 )^{1/2} (P^m_{n+1})' (t)  &=& P^{m+1}_{n+1}(t) - m(1- t^2 )^{-1/2} t P^m_{n+1} (t) , \label{eq:L2}\\
  (1 - t^2 )^{1/2} P^m_{n+1} (t)  &=& \frac{1}{2n+3}(P^{m+1}_{n+2} (t)-P^{m+1}_{n} (t)) , \label{eq:L3}\\
2 m t P^m_{n+1} (t) &=&  (1 - t^2 )^{1/2}(P^{m+1}_{n+1}(t)+(n+m+1)(n-m+2)P^{m-1}_{n+1}(t)),\label{eq:L4} \\
\label{eq:L5}
(n -m +1)P^m_{n+1} (t) &=& (2n + 1) t P^m_n (t) - (n + m)P^m_{n-1} (t)  .
\end{eqnarray}
The spherical harmonics (\ref{eq:defUV}) are polynomials when
expressed in cartesian coordinates $(x_0,x_1,x_2)$. One obtains a
basis for the space $\M^{(n)}(\B^3)$ of homogeneous monogenic
functions of degree $n$ formed by the $2n+3$ polynomials
\begin{eqnarray}
 X_m^n &=& D[\widehat{U}^{n+1}_m], \qquad m=0,\cdots,n+1 \nonumber\\
 Y_m^n &=& D[\widehat{V}^{n+1}_m], \qquad m=1,\cdots,n+1, \label{eq:XY0}
\end{eqnarray}
 These  are monogenic by construction due to the factorization of the
Laplacian.  A detailed explanation of the analogous construction for
$\R^3 \rightarrow \H$ is found in \cite{Ca,GHS}; the specific construction for $\R^3 \rightarrow \R^3$ given
here appears in \cite{Mor,MG} and it is shown that $X_m^n,Y_m^n$ form
an orthogonal basis for $\M^{(n)}(\B^3)$.  For reference we recall that
the derivation (\ref{eq:XY0}) produces the monogenic basis elements in
the following form,
\begin{eqnarray} \label{eq:XY}
  X_m^n &=& r^n\bigg( A^n_m\cos m\varphi +
  (B^n_m\cos\varphi\cos m\varphi-C^n_m\sin\varphi\sin m\varphi )e_1   \nonumber\\
  &&\ +\ (B^n_m\sin\varphi\cos m\varphi+C^n_m\cos\varphi\sin m\varphi)e_2 \bigg), \nonumber \\
  Y_m^n &=&  r^n \bigg( A^n_m\sin m\varphi +
    (B^n_m\cos\varphi\sin m\varphi+C^n_m\sin\varphi\cos m\varphi )e_1  \nonumber\\
  &&\ +\ (B^n_m\sin\varphi\sin m\varphi-C^n_m\cos\varphi\sin m\varphi)e_2 \bigg), 
\end{eqnarray}
where 
\begin{eqnarray} \label{eq:ABC}
  A^n_m &=& \frac{1}{2}\bigg( (1-t^2)(P_{n+1}^m)'(t) + (n+1)tP_{n+1}^m(t)\bigg)
    |_{t=\cos\theta}, \nonumber\\
  B^n_m &=& \frac{1}{2}\bigg( \sqrt{1-t^2} \ t(P_{n+1}^m)'(t) - (n+1) \sqrt{1-t^2}P_{n+1}^m(t) \bigg)
    |_{t=\cos\theta}, \nonumber\\
 C^n_m &=& \frac{m}{2 \sqrt{1-t^2}} P_{n+1}^m(t) 
    |_{t=\cos\theta}.
\end{eqnarray}
(We observe that other authors have used notation such as
$X_m^{n,\dagger}$, $Y_m^{n,\dagger}$ where we write $X_m^n$, $Y_m^n$, in order
to stress that these functions are defined in $\B^3$ rather than on $S^2$.) Finally in
$L_2(\B^3,\R^3)$ we consider the norm $ \|f\|_2=\sqrt{ \scpr{f}{f}}$,
where the scalar product is over $\B^3$ as in section
\ref{subsec:antiambi} (rather than over $S^2$ as is the case of some
other authors). The norms of the solid spherical harmonics and
orthogonal monogenic functions are given by
\begin{eqnarray} \label{eq:normaUV}
 \| \widehat{U}^n_0\|_2 &=&  \sqrt{\frac{4 \pi}{(2n+1)(2n+3)} } \nonumber\\
 \| \widehat{U}^n_m\|_2 &=& \| \widehat{V}^n_m\|_2 =  \sqrt{\frac{2 \pi}{(2n+1)(2n+3)} \frac{(n+m)!}{(n-m)!}},
\end{eqnarray}
\begin{eqnarray} \label{eq:XYnorms}
  \|X_0^n\|_2 &=&  \sqrt{\frac{\pi(n+1)}{2n+3}}, \nonumber\\
  \|X_m^n\|_2  &=& \|Y_m^n\|_2 \ =\  
    \sqrt{\frac{\pi(n+1)(n+m+1)!}{2(2n+3)(n-m+1)!}} ,
\end{eqnarray}
when $m\geq 1$.
 
The following explicit representation of the basis elements of
$\M^{(n)}(\B^3)$ in terms of spherical harmonics is stated in
\cite{MG} without proof.  Since our results depend on this
representation, we will prove it here in detail.

\begin{theo} \label{th:XYdeUV}
Write 
\[ c^n_m= \frac{(n+m)(n+m+1)}{4}.
\] For each degree $n\ge1$, the basis elements for the homogeneous 
monogenic polynomials of degree $n$ are given by
\begin{eqnarray} \label{eq:XYdeUV}
  X_0^n  &=& \frac{(n+1)}{2} \widehat{U}^n_0 
       + \frac{1}{2} \widehat{U}^n_1e_1  + \frac{1}{2} \widehat{V}^n_1e_2, \nonumber\\
  X_m^n &=&  \frac{(n+m+1)}{2} \widehat{U}^n_m -   \left(
     c^n_m \widehat{U}^n_{m-1} - \frac{1}{4} \widehat{U}^n_{m+1} \right)e_1
  +     \left( c^n_m \widehat{V}^n_{m-1} + \frac{1}{4}\widehat{V}^n_{m+1}\right)e_2,
     \nonumber\\
  Y_m^n &=&  \frac{(n+m+1)}{2} \widehat{V}^n_m - \left(
     c^n_m \widehat{V}^n_{m-1} - \frac{1}{4}\widehat{V}^n_{m+1}\right)e_1
  -  \left(  c^n_m \widehat{U}^n_{m-1} + \frac{1}{4}\widehat{U}^n_{m+1}\right)e_2,
\end{eqnarray}
where $1\le m\le n+1$. 
\end{theo}
\proof Recalling (\ref{eq:par0}),(\ref{eq:parvec}) we see that since
$\widehat{U}^{n+1}_m$ and $\widehat{V}^{n+1}_m$ are scalar valued, the
definition (\ref{eq:XY0}) may be expressed as
\begin{eqnarray*}
   X_m^n &=& \partial_0 \widehat{U}^{n+1}_m -
     \partial_1 \widehat{U}^{n+1}_me_1 - 
     \partial_2 \widehat{U}^{n+1}_me_2, \\
 Y_m^n &=& \partial_0 \widehat{V}^{n+1}_m -
     \partial_1 \widehat{V}^{n+1}_me_1 - 
     \partial_2 \widehat{V}^{n+1}_me_2, 
\end{eqnarray*}
so the components of $X_m^n,Y_m^n$ given in (\ref{eq:XY}) are
precisely the partial derivatives of
$\widehat{U}_m^{n+1},\widehat{V}_m^{n+1}$.  Consider the formula
proposed for $X_0^n$ in (\ref{eq:XYdeUV}); i.e., assume for the
moment $m=0$. By (\ref{eq:XY}), it is necessary to prove the three
equalities
\begin{eqnarray}
  \frac{(n+1)}{2} U^n_0   &=&  A^n_0,  \label{eq:X00}\\
   \frac{1}{2} U^n_{1} &=&  B^n_0\cos\varphi,  \label{eq:X01} \\
   \frac{1}{2} V^n_{1} &=&  B^n_0\sin\varphi.  \label{eq:X02}
\end{eqnarray}
Substitute the definitions (\ref{eq:defUV}) and (\ref{eq:ABC}) 
into these equations.  We see that  (\ref{eq:X00}) 
is given immediately by (\ref{eq:L1}). To verify equations (\ref{eq:X01}), (\ref{eq:X02}) we must prove that
\[  (1-t^2)^{1/2}t(P_{n+1})'  - (n+1)(1-t^2)^{1/2} P_{n+1}  - P_n^1  = 0.
\]
where for brevity we write $P_{n+1}$ in place of $P_{n+1}(t)$. The left
hand side may be expressed as
\[ (1-t^2)^{1/2} ( t(P_{n+1})'-(n+1) P_{n+1}- (P_n)' )
\]
and when we substitute the values for $(P_n)'$, $(P_{n+1})'$ given by  
(\ref{eq:L1}) we the result is zero according to (\ref{eq:L5}).
Thus (\ref{eq:X01}) and (\ref{eq:X02}) hold.
 
Similarly the desired formula for $X_m^n$, $1\le m\le n+1$, is equivalent to proving
the three equalities
\begin{eqnarray}
  \frac{(n+m+1)}{2} U^n_m   &=&  A^n_m\cos m\varphi , \label{eq:Xm0}\\
  -c^n_m U^n_{m-1} + \frac{1}{4} U^n_{m+1} &=&  B^n_m\cos\varphi\cos m\varphi-C^n_m\sin\varphi\sin m\varphi, \label{eq:Xm1}\\
   c^n_m V^n_{m-1} + \frac{1}{4} V^n_{m+1} &=&  B^n_m\sin\varphi\cos m\varphi+C^n_m\cos\varphi\sin m\varphi. \label{eq:Xm2}
\end{eqnarray} 
As before, substituting the definitions (\ref{eq:defUV}) and
(\ref{eq:ABC}) we obtain again that (\ref{eq:Xm0}) reduces immediately
to (\ref{eq:L1}). For equation (\ref{eq:Xm1}) we equate the
coefficients of $\cos\varphi\cos m\varphi$ and of $\sin\varphi\sin
m\varphi$ on both sides; for (\ref{eq:Xm2}) we use the coeficients of
$\sin\varphi\cos m\varphi$ and $\cos\varphi\sin m\varphi$. From this
it is seen that proving equations (\ref{eq:Xm1}) and (\ref{eq:Xm2})
is equivalent to proving
\begin{eqnarray}
 (n+m)(n+m+1)P_n^{m-1}-P_n^{m+1}+2(1-t^2)^{1/2}t(P_{n+1}^m)'  \nonumber \\
  -\ 2(n+1)(1-t^2)^{1/2}P_{n+1}^m &=& 0,   \label{eq:X1}  \\
 (n+m)(n+m+1)P_n^{m-1}+ P_n^{m+1}-2m(1-t^2)^{-1/2}P_{n+1}^m &=& 0.   \label{eq:X2}
\end{eqnarray}
Let us verify (\ref{eq:X1}). Eliminate the derivative $(P_{n+1}^m)'$
 via (\ref{eq:L1}) to convert the left hand side into
\begin{eqnarray*}
  &&\hspace{-6ex} (n+m)(n+m+1)P_n^{m-1} -
 P_n^{m+1}+2(n+m+1)t(1-t^2)^{-1/2}P_n^m  
 \\ &-& \!\!\!  2(n+1)(1-t^2)^{1/2}P_{n+1}^m  .
\end{eqnarray*}  
We note that (\ref{eq:L4}) provides a value for $(n+m)(n-m+1)P_n^{m-1}$,
so we are led to decompose $n+m+1=2(n+1)-(n-m+1)$ and to arrange the
terms as follows,
\begin{eqnarray*}
&& \hspace{-4ex} 2(n+m)(n+1)P_n^{m-1}-(P_n^{m+1}+(n+m)(n-m+1)P_n^{m-1})  \\
&&\hspace{-2ex} + \ 2(1-t^2)^{-1/2}((n+m+1)tP_n^m-(n+1)P^m_{n+1} ), \\
 &=&  2(n+m)(n+1)P_n^{m-1} - (1-t^2)^{-1/2}(
  (n-m+1)t P_n^m -2(n+1)P_{n+1}^m).
\end{eqnarray*}
Multiply this by $(1-t^2)^{1/2}$, and then substitute in the first term
the value for $(1-t^2)^{1/2}P_n^{m-1}$ provided by (\ref{eq:L3})  to arrive at
\[  (n+m)(P^m_{n+1} -P^m_{n-1}) + (2n+1)tP_n^m - (2n+1)P^m_{n+1} 
\]
which clearly vanishes by (\ref{eq:L5}), thus proving (\ref{eq:X1}).

 Now we verify equation (\ref{eq:X2}). Its left hand side can be written
as 
\begin{eqnarray*}
 2m(n+m)P_n^{m-1}-2m(1-t^2)^{-1/2}P_{n+1}^m+P_n^{m+1}+(n+m)(n-m+1)P_n^{m-1}  
\end{eqnarray*}
and by applying (\ref{eq:L4}) to the last two terms it is equal to
\[   2m(n+m)P_n^{m-1} - 2m(1-t^2)^{-1/2}(P_{n+1}^m- tP_n^m ).
\] 
We may divide by $2m(1-t^2)^{1/2}$. Using (\ref{eq:L3}) (with $n-1$,
$m-1$ in place of $m$, $n$) to replace the first term, the result is
zero by (\ref{eq:L5}). This proves (\ref{eq:X2}) and thus the formula
for $X_m^n$.  The formula for $Y_m^n$ likewise reduces to
(\ref{eq:X1}) and (\ref{eq:X2}), so (\ref{eq:XYdeUV}) is verified for
all cases.  \qed

Theorem \ref{th:XYdeUV} immediately yields  a basis for $\Vec\M$,
which turns out to be orthogonal:

\begin{prop} \label{prop:basevecM}
For each $n\ge0$, the set
\[ \{  \Vec X_m^n= \vec\partial \widehat{U}_m^{n+1},\ 0\le m\le n+1\} \cup 
   \{ \Vec Y_m^n= \vec\partial \widehat{V}_m^{n+1},\ 1\le m\le n+1\}
\]
is an orthogonal basis for $\Vec\M^{(n)}(\B^3)$. The union of these
sets over all $n\ge0$ is an orthogonal basis for $\Vec\M(\B^3)$, and
the norms of the basis elements are given by
\begin{eqnarray*}
  \|\Vec X_0^n\|  &=& \sqrt{ \frac{\pi n(n+1)}{(2n+1)(2n+3)} }\ , \\
 \|\Vec X_m^n\| &=&  \|\Vec Y_m^n\|^2 \ = \  
 \sqrt{ \frac{\pi(n^2+m^2+n)(n+m+1)!}{2(2n+1)(2n+3)(n-m+1)!} }\ , \quad m\ge1.
\end{eqnarray*}
\end{prop}

\proof The given set is clearly a basis of $\Vec\M^{(n)}(\B^3)$
because of (\ref{eq:parvec}). We need to see that it is orthogonal.
Choose elements $f,g$ in the basis for the monogenics $\{
X_0^n,\ X_m^n,\ Y_m^n\}$, and express them as $f=f_0+\vec f$,
$g=g_0+\vec g$, with scalar parts $f_0$, $g_0$.  Clearly
\[ \langle f,g\rangle = \langle f_0,g_0\rangle 
 +  \langle \vec f, \vec g\rangle.
\] 
By Theorem \ref{th:XYdeUV}, the scalar parts of the monogenics run
through (scalar multiples of) the spherical harmonics, and thus are
orthogonal.  Suppose $f\neq g$.  Then $\langle f,g\rangle=0$ and
$\langle f_0,g_0\rangle=0$, so $\langle \vec f, \vec g\rangle=0$ as
desired. To calculate the norms, now suppose $f=g$.  Then $\langle
\vec f,\vec f\rangle=\|f\|^2- \|f_0 \|^2$.  Thus if
$f=X_0^n$, then by (\ref{eq:normaUV}), (\ref{eq:XYnorms}), and
(\ref{eq:XYdeUV})
\[
 \langle \vec f,\vec f\rangle=\frac{\pi(n+1)}{(2n+3)}
 -\frac{(n+1)^2}{4} \frac{4 \pi}{(2n+1)(2n+3)}=
\frac{\pi n(n+1)}{(2n+1)(2n+3)}.
\]
Similarly, if $f=X_m^n$ or $f=Y_m^n$ for $m\ge1$,
\begin{eqnarray*}
 \langle \vec f,\vec f\rangle &=& \frac{\pi(n+1)(n+m+1)!}{2(2n+3)(n-m+1)!}
 -\frac{(n+m+1)^2}{4} \frac{2 \pi(n+m)!}{(2n+1)(2n+3)(n-m)!} \\
 &=& \frac{\pi(n^2+m^2+n)(n+m+1)!}{2(2n+1)(2n+3)(n-m+1)!}. \qed
\end{eqnarray*}

We now give a basis for the ambigenic functions. It must be noted that
the monogenic constants $X_{n+1}^n$ and $Y_{n+1}^n$ are the negatives
of their conjugates, so care must be taken to count them only
once. Thus the following functions thus form a (not orthogonal) basis
for the ambigenic functions on $\B^3$:
\begin{equation}
 X_0^n,\  X_1^n,\ \dots,\ X_n^n,\  X_{n+1}^n,\ 
     Y_1^n,\ \dots,\ Y_n^n,\  \cg{X_0^n},\
     \cg{X_1^n},\ \dots,\ \cg{X_n^n},   \  
    \cg{Y_1^n},\ \dots,\ \cg{Y_n^n},\   Y_{n+1}^n.
\end{equation}

\begin{prop} \label{prop:baseambi}
Let $n>0$. The $4n+4$ functions
\begin{eqnarray}
X^{n,+}_m &:=& X^n_m, \qquad m=0,\dots,n+1, \nonumber\\
Y^{n,+}_m &:=& Y^n_m, \qquad m=1,\dots,n, \nonumber\\
X^{n,-}_m &:=& \cg{X_m^n} - a^n_m X^n_m, \quad m=0,\dots,n, \nonumber\\
Y^{n,-}_m &:=& \cg{Y_m^n} - a^n_m Y^n_m, \quad m=1,\dots,n+1,
\end{eqnarray}
where 
\[ a^n_m=\frac{n-2m^2+1}{(n+1)(2n+1)}\quad (0\le m\le n), \quad 
   a^n_{n+1}=0,
\] 
form an orthogonal basis for the space of square integrable ambigenic
functions on $\B^3$ which are homogeneous of degree $n$.
\end{prop}

\proof Since these are $4n+4$ ambigenic functions, it suffices to
prove the orthogonality.  First we calculate the scalar products of
each $X^n_m$, $Y^n_m$ and its conjugate. By Theorem
\ref{th:XYdeUV} and (\ref{eq:normaUV}) we obtain
\begin{eqnarray} \label{eq:prodconj}
\langle X_0^n,\cg{X_0^n} \rangle &=& \frac{\pi(n+1)}{(2n+1)(2n+3)}, \nonumber\\
\langle X_m^n,\cg{X_m^n} \rangle &=& \langle Y_m^n,\cg{Y_m^n} \rangle = 
\frac{\pi(n-2m^2+1)(n+m+1)!}{2(2n+1)(2n+3)(n-m+1)!}.
\end{eqnarray}
for $1\leq m\leq n+1$.
Since the set $\{X^n_0,X^n_m,Y^n_m\colon\ m=1,\dots,n+1 \}$
is an orthogonal basis of $\M^{(n)}(\B^3)$, it follows at once that
\[
\langle X^{n,+}_m , Y^{n,+}_l \rangle =\langle X^{n,+}_m , Y^{n,-}_l \rangle =
\langle Y^{n,+}_m , X^{n,-}_l \rangle =\langle X^{n,-}_m , Y^{n,-}_l \rangle =0.
\]
Further, by (\ref{eq:XYnorms})
\[ \langle X^{n,+}_m , X^{n,+}_l \rangle = 
  \left \{ \begin{array}{ll}
    \frac{\DS \pi(n+1)}{\DS 2n+3} \qquad &\mbox{if } m=l=0, \\[1ex]
    \frac{\DS \pi(n+1)(n+m+1)!}{\DS 2(2n+3)(n-m+1)!} & \mbox{if } m=l\neq0, \\
    \quad 0 &  \mbox{if } m \neq l. \\
\end{array} \right. 
\]
Substituting (\ref{eq:XYnorms}) and (\ref{eq:prodconj}) we find that
\[   \langle X^{n,+}_m , X^{n,-}_l \rangle =\langle X^n_m ,
     \cg{X_l^n} \rangle - a^n_l\langle X^n_m , X^n_l \rangle=0 
\]
for all $l,m$, and also that
\begin{eqnarray*}
 \langle X^{n,-}_m , X^{n,-}_l \rangle &=& 
 \langle \cg{X_m^n} , \cg{X_l^n} \rangle - 
  a^n_l \langle \cg{X_m^n} , X^n_l \rangle - 
  a^n_m\langle X_m^n , \cg{X_l^n} \rangle +
  a^n_m a^n_l \langle X_m^n , X_l^n \rangle \\ &=&
  \left \{ \begin{array}{ll}
    \frac{\DS 4 \pi n(n+1)^2}{\DS (2n+3)(2n+1)^2} \qquad 
      &\mbox{if } m=l=0, \\[2ex]
    \frac{\DS 2\pi (n^2+m^2+n)(n+m+1)(n+m+1)!}
         {\DS (n+1)(2n+3)(2n+1)^2(n-m)!} & \mbox{if } m=l\neq0, \\[2ex]
\quad 0 &  \mbox{if } m \neq l. \\
\end{array} \right.
\end{eqnarray*}
The calculation of the scalar products for $\{Y_m^n\}$, as well as for the mixed
cases, is   similar.
\qed
 
\section{Contragenic functions in $\R^3$\label{sec:contragenic}}
 
\subsection{Existence of contragenics}

From now on we will abbreviate $\M^{(n)}=\M^{(n)}(\B^3)$,
$\cg\M^{(n)}=\cg\M^{(n)}(\B^3)$, $\har^{(n)}_\R=\har^{(n)}_\R(\B^3)$.
Because of the correspondence of monogenic constants with holomorphic
functions described in section \ref{sec:monofunc}, and since the real
homogeneous harmonic polynomials of degree $n$ in the complex variable
$x+iy$ are linear combinations of $\re(x+iy)^n$ and $\im(x+iy)^n$, the
space $\M^{(n)}\cap\cg\M^{(n)}$ is 2-dimensional over $\R$ for
$n\ge1$. We summarize in Table 1 the dimensions over $\R$ of the
relevant spaces of functions in $\B^3$. Recall that $\har$ denotes
$\R^3$-valued functions.
\begin{table}[ht!] \label{tab:dim}
\[
\begin{array}{|c|c|c|} \hline
 \mbox{Space of polynomials} & n=0 & \ n\ge 1 \ \\\hline
 \str \har^{(n)}_\R & 1 &2n+1 \\\hline
 \str \M^{(n)} ,\ \cg \M^{(n)} & 3 & 2n+3 \\\hline
 \str \M^{(n)} \cap \cg \M^{(n)} & 3 &2  \\\hline
 \str \M^{(n)} + \cg \M^{(n)} &  3 &4n+4 \\\hline
 \str \har^{(n)}  & 3 &6n+3 \\\hline
\end{array}
\]
\caption{Dimensions over $\R$ of spaces related to monogenic polynomials
in $\B^3$.}
\end{table}

From the last two rows of this table comes the following notable fact.
\begin{theo} \label{theo:contragenic}
Not all $L_2$-harmonic functions are ambigenic: $\M_2+\cg\M_2$ \textit{is a
  proper subspace of} $\har(\B^3)\cap L_2(\B^3)$.
\end{theo}
 
\defi In any domain $\Omega$, a harmonic function $h\in
\har(\Omega)\cap L_2(\Omega)$ is called \textit{contragenic} when it
is orthogonal to all square-integrable ambigenic functions, i.e., if
it lies in
\[  \N(\Omega) =  \big(\M_2(\Omega) + \cg \M_2(\Omega)\big)^\perp  
 \]
 where the orthogonal complement is taken in $\har(\Omega)\cap
 L_2(\Omega)$. In $\B^3$ we have the orthogonal complements in
the spaces of homogeneous harmonic polynomials $\har^{(n)}$,
\[  \N^{(n)}  = \big(\M_2^{(n)}  + \cg \M_2^{(n)} \big)^\perp.
\]

As observed in the introduction, there is no direct analogy to be
found for contragenic functions in hypercomplex analysis on $\C$ or
$\H$ or more general Clifford algebras since in those contexts all
harmonic functions are known to be ambigenic.
 
Let $n\ge1$. It follows from (\ref{eq:constmono}) that
$\M^{(n)}(\Omega)\cap\cg\M^{(n)}(\Omega)\subseteq\Vec\M^{(n)}(\Omega)$.
By Lemma \ref{lemm:f0}, when $\Omega$ is simply connected we have
$\sc\M^{(n)}(\Omega)=\har^{(n)}_\R(\Omega)$.  Returning
to $\B^3$, we have specifically $\dim\sc\M^{(n)}=\dim \har^{(n)}_\R
=2n+1$.  From (\ref{eq:dirsum}) it follows that $\dim
\Vec\M^{(n)}=2n+3$. (Since this is equal $\dim\M^{(n)}$, this means
that given the vectorial part $\vec f$, the scalar part $f_0$ is
uniquely determined.)  Thus
\begin{equation} \dim\N^{(n)} = 2n-1 
\end{equation}
when $n\ge1$, while $\dim\N^{(0)} = 0$.

The following result is a simple consequence.
\begin{theo}
Let $h\in\N=\N(\B^3)$. Then $h$ can be uniquely expressed as a sum in
$L_2$
\[   h = \sum_{n=1}^\infty h^{(n)}
\]
where $h^{(n)}\in\N^{(n)}$.
\end{theo}

\note It is easily checked by a dimension count that the analogue of Theorem
\ref{theo:contragenic} for ``clasical" monogenic functions $\H
\rightarrow \H$ does not hold; i.e., all harmonics are ambigenic over $\H$. For
$n\geq0$ one has that the homogeneous monogenics in $\B^4$ form a
right vector space over $\H$ of dimension $\frac{1}{2}(n+1)(n+2)$
\cite{Sud}, and the same is true for the antimonogenics. The monogenic
constants have dimension $n+1$ over $\H$ (See example \cite{Bock2}).
Since the dimension of the harmonics from $\H$ to $\H$ is
$(n+1)^2$ over $\H$ \cite{Sud}, it follows that every harmonic
function can be expressed as a sum of a monogenic function and an
antimonogenic function.

It may also shed light on the situation to see what fails when one
attempts to express a harmonic $\R^3$-valued function in terms of
monogenics and antimonogenics.  For scalar-valued
$f_0\in\har_\R(\Omega)$ there is in fact no problem, since
\[ f_0 = \frac{1}{2}(f_0+\vec g_0) +  \frac{1}{2}(f_0-\vec g_0)
\]
where $f_0+\vec g_0$ is the completion of $f_0$ to a monogenic function as
given by Lemma \ref{lemm:f0}. When we are given a general harmonic
function $f_0+f_1e_1+f_2e_2\in\har(\Omega)$ and complete each component separately,
we obtain analogously
\begin{eqnarray*}
 f_0+f_1e_1+f_2e_2 &=&
  \frac{1}{2}( (f_0+\vec g_0)+(f_1+\vec g_1)e_1+(f_2+\vec g_2)e_2 ) \\
 && +   \frac{1}{2}( (f_0-\vec g_0)+(f_1-\vec g_1)e_1+(f_2-\vec g_2)e_2 ). 
\end{eqnarray*}
However, the two functions added on the right hand side of this equation, while in 
the kernels of $\Dbar$ and $D$ respectively,
need not take their values in $\R^3$.  Thus we do not obtain in this
natural way a representation as an ambigenic function.

\subsection{Basic properties of contragenics}

Although our main interest is in $\B^3$,
we observe also some relations which hold in more general domains.
Consider an arbitrary contragenic function $h= h_0e_0 + h_1e_1 +
h_2e_2\in\N(\Omega)$.  By the orthogonal decomposition (\ref{eq:dirsum}),
$h\perp\sc\M_2(\Omega)$, which implies $\scpr{h_0}{h_0}=0$ so in fact 
$h_0=0$. Thus $h$ is of the form
\begin{equation}
 h(x_0,x_1,x_2)=   h_1(x_0,x_1,x_2)e_1 + h_2(x_0,x_1,x_2)e_2 .
\end{equation}
In particular, contragenic functions are never invertible.
Equation (\ref{eq:dirsum}) also gives the property
$h\perp\Vec\M_2(\Omega)$, which by Lemma \ref{lemm:f0} can be expressed
as
\begin{equation} \label{eq:hperpf}
  \int_\Omega(h_1f_1+h_2f_2)\,dV=0 
\end{equation}
whenever $\partial_1 f_2=\partial_2 f_1$ with $f_1$ and $f_2$ harmonic.

Work on monogenic functions in different contexts has
focused on standard domains such as spheres, ellipsoids, cylinders,
rectangles, etc.\ (see for example \cite{ML}).  Let $\Omega_1\subseteq\R^3$ be a domain which
enjoys the symmetry that $(x_0,x_2,x_1)\in\Omega_1$ whenever
$x=(x_0,x_1,x_2)\in\Omega_1$. For $f:\Omega_1\to\H$ write
\[ f^*(x)=f_0(x_0,x_2,x_1)+f_2(x_0,x_2,x_1)e_1+f_1(x_0,x_2,x_1)e_2-f_3(x_0,x_2,x_1)e_3
\]
for $x\in\Omega_1$.  Thus $(f^*)^*=f$, and it is easily seen that
\begin{equation} 
    (fg)^*= g^*f^*
 \end{equation}
for $f,g:\Omega_1\rightarrow\R^3$, and
\begin{equation} \label{eq:Dbarf*}
  \Dbar(f^*) = (\Dbar f)^*
\end{equation}
when $f$ is differentiable. From this we have
 
\begin{prop}  \label{prop:flip}
 The involution $f\mapsto f^*$ preserves $\M(\Omega_1)$ and
 $\cg\M(\Omega_1)$. If $\Omega_1$ is simply connected, the involution
preserves $\N(\Omega_1)$ as well.
\end{prop}

\proof Let $f\in\M(\Omega_1)$. By (\ref{eq:Dbarf*}),
$\Dbar(f^*)=0^*=0$, so $f^*\in\M(\Omega_1)$. Similarly $\cg\M$ is
invariant. Now let $h\in \N(\Omega_1)$.  For $g\in
\har_{\R}(\Omega_1)\cap L_2(\Omega_1)$ take
$f=\partial_1g\,e_1+\partial_2g\,e_2\in \Vec\M_2(\Omega_1)$.  We find
that 
\begin{eqnarray*}
 \scpr{h^*}{f} &=& \int_{\Omega_1} \big(
     h_2(x_0,x_2,x_1)\,\partial_1 g(x)
   + h_1(x_0,x_2,x_1)\,\partial_2 g(x) \big) \,dV\\
 &=& \int_{\Omega_1} \big(
     h_2(x_0,x_2,x_1)\,\partial_2 g^*(x_0,x_2,x_1)
   + h_1(x_0,x_2,x_1)\,\partial_1 g^*(x_0,x_2,x_1) \big)\,dV
 \ =\ 0 
\end{eqnarray*}
by (\ref{eq:hperpf}), where $g^*(x)=g(x_0,x_2,x_1)$. Assuming
$\Omega_1$ is simply connected, it follows from Lemma \ref{lemm:f0}
that $f$ ranges over all of $\Vec\M_2(\Omega_1)$; thus $h^*\in\N(\Omega_1)$ as claimed.
\qed
 
We conjecture that the simple connectedness is not necessary for
the invariance of $\N(\Omega_1)$ in Proposition \ref{prop:flip}.

In the unit ball we have the following characterization of contragenic
functions via integration over the unit sphere.  

\begin{prop} \label{prop:ncriterion}
Let $h\in \har^{(n)}(\B^3)$, $h=h_1e_1+h_2e_2$. Then $h\in\N^{(n)}(\B^3)$
if and only if the equality
\begin{equation}\label{eq:S2condition}
 \int_{S^2}h_1 g\, dx_0\wedge dx_2=\int_{S^2}h_2 g\, dx_0\wedge dx_1 
\end{equation}
holds for all $g\in\har^{n+1}_\R(\B^3)$.
\end{prop}
\proof We have the following relation of  differential forms,
\[  (h_1\partial_1 g + h_2\partial_2 g)\,dx_0\,dx_1\,dx_2  
 =\hspace*{35ex}
\] \[  \hspace*{15ex}
   d(-h_1 g\, dx_0\wedge dx_2+h_2 g\, dx_0\wedge dx_1)
   -(\partial_1 h_1 + \partial_2 h_2)g\,dx_0\,dx_1\,dx_2.
\]
Since harmonic functions of differing degrees are orthogonal over
$\B^3$, integration leaves
\[  \scpr{h}{\vec\partial g} = \int_{\B^3} d(-h_1 g \,dx_0\wedge dx_2+h_2 g \,dx_0\wedge dx_1)  
\]
 which is equal to 0 whenever $h\in\N$, and
Stokes' theorem gives (\ref{eq:S2condition}). Conversely, if
(\ref{eq:S2condition}) holds, then $\langle h, \vec\partial
g\rangle=0$. Since $\vec\partial g$ ranges over all of $\Vec\M^{(n)}$, $h$
is orthogonal to $\M^{(n)}$.  Further, since $h$ is trivially orthogonal to
$\M^{(m)}$ for $m\not=n$, it is contragenic.  \qed
    
\begin{coro} \label{coro:hcriterion}
Let $h\in \har(\cg\B^3)$, $h=h_1e_1+h_2e_2$. Then $h\in\N(\B^3)$ if and only if (\ref{eq:S2condition}) holds  for all $g\in
\har_\R(\cg\B^3)$.
\end{coro} 

Proposition \ref{prop:ncriterion} and Corollary \ref{coro:hcriterion},
together with the bases given in section \ref{sec:homomono},
provide an algorithmic method for determining when a given harmonic
function, expressed as a convergent series in $\H\cap L_2$, is
contragenic.  However, it is not likely that there is a simple
characterization of contragenics purely in terms of derivatives for
general $\Omega$, and indeed we know of none even for $\B^3$.
   
\subsection{Bergman kernel for $\Vec\M$} 

The natural generalization of the holomorphic Bergman kernel \cite{Hil}
from the context of holomorphic functions in $\C$ to that of monogenic functions
in $\H$ is described in \cite{SV1,SV2}.  A generalization for functions
in $\R^3$ was defined and studied more recently in \cite{GLS}.  We
restate some of the main facts in the present terminology, and then
give a new ``Bergman kernel'' which is more appropriate to the
subject at hand, as it provides another characterization of contragenic
functions.  

The following result establishes that evaluation at a fixed
point is a continuous linear functional on $\Vec\M(\Omega)$; it suffices to
work in $\B^3$ for the basic estimate.

\begin{prop} \label{prop:evalcont}
Let $f\in\Vec\M$.  Then for all $x\in\B^3$,
\[ |f(0)| \le C\|f\|_2 
\]
where $C=\sqrt{3/(4\pi)}$.
\end{prop}

\proof Since $f$ is harmonic, using the orthogonal basis of
Proposition \ref{prop:basevecM} we see that the constant $f(0)$ is
orthogonal to $f-f(0)$ in $L_2(\B^3)$.  (More simply, one may just
observe that the constants are orthogonal to all harmonic functions
which vanish at the origin, a statement of the Mean Value Property of
harmonic functions.) Therefore
\[  |f(0)|^ 2 = C^2\int_{B^3}|f(0)|^2\,dV =
     C^2 ( \|f\|_2^2 -\|f-f(0)\|_2^2  ) \le C^2\|f\|_2^2 . \qed
\]

The underlying idea, given a closed subspace $\A\subseteq L_2$ and an
orthonormal basis $\{\varphi_k\}$ of $\A$, is to form an integral kernel
$ B(x,y)=\sum_k \varphi_k(x) \cg{\varphi_k(y)}$ which automatically enjoys the reproducing property
$  f(x) = \int B(x,y)\,\cg{f(y)}\,dy $
for $f\in\A$, and projects $L_2$ orthogonally onto $\A$.
 However, given an orthornormal basis $\{\varphi_k\}$ of
$\M_2(\Omega)$, there is a problem if we try to construct a Bergman
kernel in this way because $\M_2$ is not closed under multiplication:
the integrand of
\begin{eqnarray*}
     \int \big( B_0(x,y)e_0+B_1(x,y)e_1+B_2(x,y)e_2+B_3(x,y)e_3) 
   \cdot (f_0(y)e_0+f_1(y)e_1+f_2(y)e_2 \big)\,dV_y 
\end{eqnarray*}
contains the term $(B_1f_2-B_2f_1+B_3f_0)e_3$, whereas $f$ should be 
$\R^3$-valued.  Thus one needs an additional condition
\begin{equation}\label{eq:intBf} 
  \int(B_1f_2-B_2f_1+B_3f_0)\,dy=0. 
\end{equation}
In \cite{GLS} this is dealt with by working in the subspace of $L_2$
corresponding to functions $f$ for which this property holds.  We will
take a different approach here, constructing a Bergman kernel for
$\Vec\M_2$ rather than $\M_2$, and not requiring a special
condition such as (\ref{eq:intBf}). 

The scalar product restricted to
$L_2(\Omega,\{0\}\oplus\R^2\oplus\{0\})$, a Hilbert space which contains
$\Vec\M_2(\Omega)$ as a closed subspace, is
\[ \scpr{f}{g}=\int(f_1g_1+f_2g_2)=-\sc\int fg 
\]
since $\cg f=-f$.  Take an orthonormal basis $\{\psi_k\}_{k=1}^\infty$
of $\Vec\M_2(\Omega)$ over $\R$, and write
$\psi_k=\psi_{k,1}e_1+\psi_{k,2}e_2$.  Define the following functions
$\Omega\times\Omega\to\{0\}\oplus\R^2\oplus\{0\}$,
\begin{eqnarray} 
  b_{\Omega,1}(x,y) &=& -\sum_{k=1}^\infty \psi_{k,1}(x) \psi_k(y), \nonumber\\ 
  b_{\Omega,2}(x,y) &=& -\sum_{k=1}^\infty \psi_{k,2}(x) \psi_k(y);\label{eq:b12}
\end{eqnarray}
it can shown by means of Proposition \ref{prop:evalcont} that these
series converge uniformly on compact subsets of $\Omega\times\Omega$.

\defi \label{def:bergmanop}
The \textit{Bergman operator $B_\Omega$ for $\Vec\M(\Omega)$} is defined by
\[
    B_\Omega[f](x) = \sc\left(\int_\Omega b_{\Omega,1}(x,y)f(y)\,dV_y\right)e_1
           +\sc\left(\int_\Omega b_{\Omega,2}(x,y)f(y)\,dV_y\right)e_2
\]  
for all $f\in L_2(\Omega,\{0\}\oplus\R^2\oplus\{0\})$ and $x\in\Omega$.

It is shown in the traditional way that $b_{\Omega,1}(x,y)$ and
$b_{\Omega,2}(x,y)$ are independent of the orthonormal basis chosen.
Since we can express elements of $\Vec\M$ as $f=\sum a_k\psi_k$
($a_k\in\R$), the following reproducing property is easily checked.
\begin{theo} \label{theo:bergman} The 
linear operator $B_\Omega$ projects
$L_2(\Omega,\{0\}\oplus\R^2\oplus\{0\})$ orthogonally onto
$\Vec\M_2(\Omega)$.  In particular, $B_\Omega[f]=f$ for
$f\in\Vec\M_2(\Omega)$.
\end{theo}
 
Note that two separate integral kernels (\ref{eq:b12}) are necessary in the
definition of the operator $B_\Omega$ because the scalar product on
$\Vec\M(\Omega)$ is only bilinear over the reals.
 
For $\B^3$, in terms of the specific basis $u^n_m :=\vec\partial
\widehat{U}^n_m$, $v^n_m :=\vec\partial \widehat{V}^n_m$ of $\Vec\M^{(n)}$ given in
Proposition \ref{prop:basevecM}, one can define kernels for each degree $n$,
\begin{eqnarray} \label{eq:bn}
  b^n_1(x,y) &=& - \sum_{m=0}^n u^n_{m,1}(x) u^n_m(y)
  - \sum_{m=1}^n  v^n_{m,1}(x) v^n_m(y), \nonumber\\ 
  b^n_2(x,y) &=& - \sum_{m=0}^n u^n_{m,2}(x) u^n_m(y)
  - \sum_{m=1}^n  v^n_{m,2}(x) v^n_m(y),  
\end{eqnarray}
and then can form operators $B^{(n)}$  analogously to Definition
\ref{def:bergmanop}.  These operators project the
harmonic functions onto $\Vec\M^{(n)}$, and the Bergman operator $B$
for $\B^3$ is their sum $B=\sum_{n=0}^\infty B^{(n)}$.
It would be interesting to express (\ref{eq:bn}) in closed form.

The following is an immediate consequence of the foregoing, and with
the formulas (\ref{eq:bn}) allows one to detect computationally when a
harmonic function is contragenic, or close to contragenic in the
$L_2$-sense.

\begin{coro} 
Let $h=h_1e_1+h_2e_2\in L_2(\Omega)$ be harmonic. Then
$h\in\N(\Omega)$ if and only if $B_\Omega[h]=0$.
\end{coro}
  
\section{Construction of homogeneous contragenic polynomials\label{sec:bases}}

In this section we will give an explicit construction of 
a basis of $\N^{(n)}$ for every $n=0,1,\dots$

One possible approach would be as follows. The basis for
$\M^{(n)}+\cg\M^{(n)}$ given in Proposition \ref{prop:baseambi} is
orthogonal, so one may extend this to a basis of $\har^{(n)}$ by
choosing suitable linearly independent triples of spherical harmonics and
applying the Gram-Schmidt process to produce the contragenic
polynomials of degree $n$. However, this procedure is quite costly
numerically,\footnote{Calculations in \textit{Mathematica} 
on a desktop computer with 4 Gb of RAM have saturated the memory when
attempting to calculate contragenic homogeneous polynomials of degree
$n\ge7$ via Gram-Schmidt as described here.} and leads to little
insight regarding contragenic functions.

Here we give a direct construction of the contragenic homogeneous
functions. From Table \ref{tab:dim} it is clear that $\N^{(0)}=\{0\}$.
  
\begin{theo} Let $n\ge1$. Write $d_m^n=(n-m)(n-m+1)$. The $2n-1$ functions
\begin{eqnarray} 
  Z^n_{0}\  \  &:=& \widehat{V}^n_1 e_1 - \widehat{U}^n_1 e_2, \nonumber \\
  Z^n_{m,+} &:=& ( d_m^n \widehat{V}^n_{m-1} + \widehat{V}^n_{m+1})e_1 +
              ( d_m^n\widehat{U}^n_{m-1} - \widehat{U}^n_{m+1})e_2  ,\nonumber  \\
  Z^n_{m,-} &:=& ( d_m^n\widehat{U}^n_{m-1} + \widehat{U}^n_{m+1})e_1 +
              (-d_m^n\widehat{V}^n_{m-1} + \widehat{V}^n_{m+1})e_2  ,  \label{eq:contragenics} 
\end{eqnarray}
 for $1\le m\le n-1$, form an orthogonal basis for $\N^{(n)}$ over $\R$.

\end{theo}

\proof  First we show that the functions (\ref{eq:contragenics}) are contragenic:
it is sufficient to show that each
one is orthogonal to $\M^{(n)}\cup\cg{\M}^{(n)}$.  As
we have already noted, since they have no scalar parts it suffices to show that
each one is orthogonal to $\Vec\M^{(n)}$, and to do
this, by (\ref{eq:parvec}) we may use the basis $\Vec\M^{(n)}$
obtained by dropping the scalar parts of the basis for $\M^{(n)}$
given by Theorem \ref{th:XYdeUV}.  Let $1\le m\le n-1$.  From Theorem
\ref{th:XYdeUV},
\[
 \langle Z^n_{m,+} ,\Vec X^n_0 \rangle \ = \
\frac{1}{2}( \langle d_m^n\widehat{V}^n_{m-1},\widehat{U}^n_1 \rangle +
\langle \widehat{V}^n_{m+1},\widehat{U}^n_1\rangle  +
\langle   d_m^n\widehat{U}^n_{m-1},\widehat{V}^n_1 \rangle -\langle \widehat{U}^n_{m+1},\widehat{V}^n_1\rangle ) .
\]
By the orthogonality of the spherical harmonics $\widehat{U}^n_m$ and
$\widehat{V}^n_m$, this scalar product is equal to zero. Next we
observe that for $1\le k\le n$,
\begin{eqnarray*}
 \langle Z^n_{m,+} ,\Vec X^n_k \rangle &=&
  \langle d_m^n\widehat{V}^n_{m-1} + \widehat{V}^n_{m+1},\ 
          -c^n_k \widehat{U}^n_{k-1} + \frac{1}{4} \widehat{U}^n_{k+1}   \rangle  \\
  && + \   \langle  d_m^n\widehat{U}^n_{m-1} - \widehat{U}^n_{m+1},\
                    c^n_k \widehat{V}^n_{k-1} + \frac{1}{4}\widehat{V}^n_{k+1} \rangle \\
 &=&  0,
\end{eqnarray*}
since $\widehat{U}^n_{m\pm1}$ is orthogonal to $\widehat{V}^n_{m\pm1}$.

Finally, it remains to check that
\begin{eqnarray*}
 \langle Z^n_{m,+} ,\Vec Y^n_k \rangle &=&
  - d_m^nc^n_k \langle \widehat{V}^n_{m-1},\widehat{V}^n_{k-1}  \rangle 
  + \frac{1}{4} d_m^n \langle \widehat{V}^n_{m-1},\widehat{V}^n_{k+1} \rangle
 \\ && - \ c^n_k \langle \widehat{V}^n_{m+1},\widehat{V}^n_{k-1}  \rangle 
  + \frac{1}{4}  \langle \widehat{V}^n_{m+1},\widehat{V}^n_{k+1}  \rangle \\
&& - \ d_m^nc^n_k \langle \widehat{U}^n_{m-1},\widehat{U}^n_{k-1}  \rangle
  - \frac{1}{4} d_m^n \langle \widehat{U}^n_{m-1},\widehat{U}^n_{k+1}  \rangle
 \\ && + \ c^n_k \langle \widehat{U}^n_{m+1},\widehat{U}^n_{k-1}  \rangle 
  + \frac{1}{4}  \langle \widehat{U}^n_{m+1},\widehat{U}^n_{k+1}  \rangle  \\
 &=&  0.
\end{eqnarray*} 
Once again it is immediate that this is true under the condition
$k\neq m$, $k\neq m+2$ and $k\neq m-2$ since all of the scalar
products involved vanish. Now consider $k=m$. Substituting the
equation (\ref{eq:normaUV}) we obtain that
\begin{eqnarray*}
 \langle Z^n_{m,+} ,\Vec Y^n_m \rangle &=&
  -d_m^nc^n_m \langle \widehat{V}^n_{m-1},\widehat{V}^n_{m-1}  \rangle
  - \ d_m^nc^n_m \langle \widehat{U}^n_{m-1},\widehat{U}^n_{m-1}  \rangle \\
&& + \frac{1}{4} \langle \widehat{V}^n_{m+1},\widehat{V}^n_{m+1} \rangle
   + \frac{1}{4} \langle \widehat{U}^n_{m+1},\widehat{U}^n_{m+1} \rangle \\
&=& -2 (n-m)(n-m+1) \frac{(n+m)(n+m+1)}{4(2n+1)(2n+3)}
    \frac{2 \pi(n+m-1)!}{(n-m+1)!} \\
&& +  \frac{2\cdot2 \pi}{4(2n+1)(2n+3)} \frac{(n+m+1)!}{(n-m-1)!} \\
&=&  0.
\end{eqnarray*}
When $k=m+2$, we see by (\ref{eq:normaUV}) that
\begin{eqnarray*}
 \langle Z^n_{m,+} ,\Vec Y^n_{m+2} \rangle &=&
  - c^n_{m+2} \langle \widehat{V}^n_{m+1},\widehat{V}^n_{m+1}  \rangle
  + c^n_{m+2} \langle \widehat{U}^n_{m+1},\widehat{U}^n_{m+1}  \rangle \\
 &=&  0
\end{eqnarray*}
and the case $k=m-2$ is similar.  Therefore $\langle Z^n_{m,+} ,\Vec
Y^n_k \rangle=0$ for all $k$.  This completes the proof that
$Z^n_{m,+}$ is contragenic.  The proofs that $Z^n_0$ and $Z^n_{m,-}$
are contragenic are analogous.

Now we show that these functions form an orthogonal basis.  Since
$\dim_\R\N^{(n)}=2n-1$, it suffices to show that they form an
orthogonal collection.  The only nontrivial cases are $\langle
Z^n_0,Z^n_{2,+} \rangle$, $\langle Z^n_{m,+},Z^n_{m+2,+} \rangle$,
$\langle Z^n_{m,+},Z^n_{m-2,+} \rangle$, $\langle
Z^n_{m,-},Z^n_{m+2,-} \rangle$ and $\langle Z^n_{m,-},Z^n_{m-2,-}
\rangle$, all of which are seen to be zero by repeated applications of
(\ref{eq:normaUV}).   \qed
 
\begin{coro}
 The set 
\[ \{ Z^n_{0},\ Z^n_{m,\pm},\ 1\le m\le n-1,\ n\ge1 \}
\] is an orthogonal basis for $\N=\N(\B^3)$.  The norms of the
basis elements are
\begin{eqnarray*}
  \| Z^n_{0} \| &=&   \sqrt{ \frac{4 \pi n(n+1)}{(2n+1)(2n+3)}}, \\
  \| Z^n_{m,\pm}\| &=&   \sqrt{\frac{8 \pi (n^2+m^2+n)(n+m-1)!}{(2n+1)(2n+3)(n-m-1)!}},
\end{eqnarray*}
for $1 \leq m \leq n-1$.
\end{coro}

\proof It is only necessary to establish the values of the norms. By (\ref{eq:normaUV}) 
\begin{eqnarray*}
\| Z^n_{0} \|^2 = \langle Z^n_0 ,Z^n_0 \rangle = \langle \widehat{V}^n_1 ,\widehat{V}^n_1 \rangle + \langle \widehat{U}^n_1 ,\widehat{U}^n_1 \rangle =
\frac{4\pi n(n+1)}{(2n+1)(2n+3)}.
\end{eqnarray*}
and,
\begin{eqnarray*}
 \| Z^n_{m,+} \|^2 &=& \langle Z^n_{m,+} ,Z^n_{m,+} \rangle \\
&=& (d^n_m)^2 \langle \widehat{V}^n_{m-1} ,\widehat{V}^n_{m-1} \rangle  + (d^n_m)^2 \langle \widehat{U}^n_{m-1} ,\widehat{U}^n_{m-1} \rangle
 + \langle \widehat{V}^n_{m+1} ,\widehat{V}^n_{m+1} \rangle + \langle \widehat{U}^n_{m+1} ,\widehat{U}^n_{m+1} \rangle \\
&=& \frac{4 \pi (n-m)^2(n-m+1)^2 (n+m-1)!}{(2n+1)(2n+3)(n-m+1)!} +  \frac{4 \pi (n+m+1)!}{(2n+1)(2n+3)(n-m-1)!} \\
&=& \frac{8 \pi (n^2+m^2+n)(n+m+1)!}{(2n+1)(2n+3)(n+m)(n+m+1)(n-m-1)!}\\
&=& \frac{8 \pi (n^2+m^2+n)(n+m-1)!}{(2n+1)(2n+3)(n-m-1)!}.
\end{eqnarray*}
The calculation for $Z^n_{m,-}$ is similar.
\qed

\note The involution  $f\mapsto f^*$ of Proposition \ref{prop:flip}
sends $Z^n_{m,+}$ to (a multiple of) $Z^n_{m,-}$ for some, but not all $m$.

\section{Conclusions}

Consider a triple $f=f_0e_0+f_1e_1+f_2e_2$ of harmonic functions in a
domain $\Omega$.  We have shown that $f$ has a natural decomposition
$f=g+h$ where $g$ is ambigenic and $h$ is orthogonal in $L_2(\Omega)$
to all ambigenic functions. The existence of nontrivial contragenic
functions raises the following question.
Suppose that $\Omega$ has smooth boundary and $f$ is defined only on
$\partial\Omega$. When is the harmonic extension of $f$ to the
interior of $\Omega$ monogenic, ambigenic, or contragenic? How do the
boundary values of the monogenic, ambigenic, or contragenic part of
the extension relate to the original $f$?

Further, it remains to investigate bases of contragenic functions in domains of
$\R^3$ other than $\B^3$, as well as analogous notions of contragenicity with respect
to other scalar products, for example in weighted inner product spaces or with respect
to the Fischer product \cite{Ros,Som}.

\end{document}